\documentclass[12pt,a4paper]{amsart}

\usepackage{anysize}
\usepackage[latin1]{inputenc}
\usepackage{hyperref}

\marginsize{3cm}{3cm}{3.5cm}{3.5cm}

\begin{document}
\title{WZ-proofs of ``divergent'' Ramanujan-type series}
\author{Jesús Guillera}
\email{jguillera@gmail.com}
\address{Av. Cesáreo Alierta, 31 esc. izda {\rm $4^o$}--A, Zaragoza (Spain)}
\date{}

\newtheorem{observacion}{Remark}

\begin{abstract}
We prove some ``divergent'' Ramanujan-type series for $1/\pi$ and $1/\pi^2$ applying a Barnes-integrals strategy of the WZ-method.
\end{abstract}

\dedicatory{Dedicated to Herb Wilf on his $80^{th}$ birthday}

\maketitle

\section{Wilf-Zeilberger's pairs}
We recall that a function $A(n,k)$ is \emph{hypergeometric} in its two variables if the quotients
\[
\frac{A(n+1,k)}{A(n,k)} \quad {\rm and} \quad \frac{A(n,k+1)}{A(n,k)}
\]
are rational functions in $n$ and $k$, respectively. Also, a pair of hypergeometric functions in its two variables, $F(n,k)$ and $G(n,k)$, is said to be a \emph{Wilf and Zeilberger (WZ) pair} \cite[Chapt. 7]{petkovsek} if
\begin{equation}\label{pro-WZ-pair}
F(n+1,k)-F(n,k)=G(n,k+1)-G(n,k).
\end{equation}
In this case, H. S. Wilf and D. Zeilberger \cite{wilf} have proved that there exists a rational function $C(n,k)$ such that
\begin{equation}\label{certificado}
G(n,k)=C(n,k)F(n,k).
\end{equation}
The rational function $C(n,k)$ is the so-called \emph{certificate} of the pair $(F,G)$. To discover WZ-pairs, we use Zeilberger's Maple package EKHAD \cite[Appendix A]{petkovsek}. If EKHAD certifies a function, we have found a WZ-pair! We will write the functions $F(n,k)$ and $G(n,k)$ using rising factorials, also called Pochhammer symbols, rather than the ordinary factorials. The rising factorial is defined by
\begin{equation}\label{poch1}
(x)_n=\left\{
\begin{array}{ll}
x(x+1)\cdots(x+n-1), & \quad n \in \mathbb{Z}^{+}, \\
1, & \quad n=0,
\end{array} \right.
\end{equation}
or more generally by $(x)_t=\Gamma(x+t)/\Gamma(x).$ For $t \in \mathbb{Z}-\mathbb{Z}^{-}$, this last definition coincide with (\ref{poch1}). But it is more general because it is also defined for all complex $x$ and $t$ such that
$x+t \in \mathbb{C}- ( \mathbb{Z}-\mathbb{Z}^{+} )$.

\section{A Barnes-integrals WZ strategy}

If we sum (\ref{pro-WZ-pair}) over all $n \geq 0$, we get
\begin{equation}\label{sumas-wz}
\sum_{n=0}^{\infty} G(n,k) - \sum_{n=0}^{\infty} G(n,k+1) = -F(0,k) + \lim_{n \to \infty} F(n,k)
\end{equation}
whenever the series above are convergent and the limit is finite. D. Zeilberger was the first to apply the WZ-method to prove a Ramanujan-type series for $1/\pi$ \cite{ekhad}. Following his idea, in a series of papers \cite{guilleraAAMwz}, \cite{guilleraRJgen}, \cite{guillera-RJonwz}, \cite{guillera-RJnewrama} and in the author's thesis \cite{guilleratesis}, we use WZ-pairs together with formula (\ref{sumas-wz}) to prove a total of eleven Ramanujan-type series for $1/\pi$ and four Ramanujan-like series for $1/\pi^2$. However, while we discovered those pairs we also found some WZ-pairs corresponding to ``divergent'' Ramanujan-type series \cite{gui-zu}, like the following pair:
\[
F(n,k)=A(n,k) \frac{(-1)^n}{\Gamma(n+1)} \left( \frac{16}{9} \right)^n, \quad
G(n,k)= B(n,k) \frac{(-1)^n}{\Gamma(n+1)} \left( \frac{16}{9} \right)^n,
\]
where
\[ A(n,k)=U(n,k)\frac{-n(n-2)}{3(n+2k+1)}, \quad B(n,k)=U(n,k)(5n+6k+1), \]
and
\[
U(n,k)=\frac{\left( \frac{1}{2} \right)_n \left( \frac{1}{4}+\frac{3k}{2} \right)_n \left( \frac{3}{4}+\frac{3k}{2} \right)_n}{(1+k)_n (1+2k)_n}  \frac{\left( \frac{1}{6} \right)_k \left( \frac{5}{6} \right)_k}{(1)_k^2}.
\]
We cannot use formula (\ref{sumas-wz}) with this pair because the series is divergent and the limit is infinite, due to the factor $(-16/9)^n$. To deal with this kind of WZ-pairs we will proceed as follows: First we replace the factor $(-1)^n$ with $\Gamma(n+1)\Gamma(-n)$. By doing it we again get a WZ-pair, because $(-1)^n$ and $\Gamma(n+1)\Gamma(-n)$ transform formally in the same way under the substitution $n \rightarrow n+1$; namely, the sign changes. To fix ideas, the modified version of the WZ-pair above is
\[
\widetilde{F}(s,t)=A(s,t) \Gamma(-s) \left( \frac{16}{9} \right)^s, \quad
\widetilde{G}(s,t)=B(s,t) \Gamma(-s) \left( \frac{16}{9} \right)^s.
\]
Then, integrating from $s=-i\infty$ to $s=i\infty$ along a path $\mathcal{P}$ (curved if necessary) which separates the poles of the form $s=0,1,2 \dots$ from all the other poles, we obtain
\begin{equation}
\frac{1}{2 \pi i} \int_{-i\infty}^{i\infty}B(s,t)\Gamma(-s)(-z)^sds=
\sum_{n=0}^{\infty}B(n,t)\frac{z^n}{n!}, \quad |z|<1,
\end{equation}
where we have used the Barnes integral theorem, which is an application of Cauchy's residues theorem using a contour which closes the path with a right side semicircle of center at the origin and infinite radius. The Barnes integral gives the analytic continuation of the series to $z \in \mathbb{C}-[1,\infty)$.
Integrating along the same path the identity $\widetilde{G}(s,t+1)-\widetilde{G}(s,t)=\widetilde{F}(s+1,t)-\widetilde{F}(s,t)$, we obtain
\begin{align}
 \int_{-i\infty}^{i\infty} \widetilde{G}(s,t+1)ds- \int_{-i\infty}^{i\infty} \widetilde{G}(s,t)ds =&
\int_{-i\infty}^{i\infty} \widetilde{F}(s+1,t)ds- \int_{-i\infty}^{i\infty} \widetilde{F}(s,t)ds  \\
= \nonumber \int_{1-i\infty}^{1+i\infty} \widetilde{F}(s,t)ds- \int_{-i\infty}^{i\infty} \widetilde{F}(s,t)ds=&-\int_{\mathcal{C}} \widetilde{F}(s,t)ds, \nonumber
\end{align}
where $\mathcal{C}$ is the contour limited by the path $\mathcal{P}$, the same path but moved one unit to the right, and the lines $y=-\infty$ and $y=+\infty$. As the only pole inside this contour is at $s=0$ and the residue at this point is zero, the last integral is zero and we have
\begin{equation}
\int_{-i\infty}^{i\infty}\widetilde{G}(s,t)ds=\int_{-i\infty}^{i\infty} \widetilde{G}(s,t+1)ds.
\end{equation}
This implies, by Weierstrass's theorem \cite{whittaker}, that
\begin{align}\nonumber
\frac{1}{2 \pi i} \int_{-i\infty}^{i\infty} \widetilde{G}(s,t)ds &=
\lim_{t \to \infty} \frac{1}{2 \pi i} \int_{-i\infty}^{i\infty} \widetilde{G}(s,t)ds=
\frac{1}{2 \pi i} \int_{-i\infty}^{i\infty} \lim_{t \to \infty} \widetilde{G}(s,t)ds \\ &=
\frac{1}{2 \pi i} \int_{-i\infty}^{i\infty} \frac{3}{\pi} \left( \frac{1}{2} \right)_{\!s} \Gamma(-s) 2^s ds=
\frac{\sqrt{3}}{\pi}, \nonumber
\end{align}
where the last equality holds because
\[
\frac{1}{2 \pi i} \int_{-i\infty}^{i\infty} \left( \frac{1}{2} \right)_{\!s} \Gamma(-s) (-z)^s ds =\sum_{n=0}^{\infty} \frac{\left( \frac12\right)_n}{(1)_n}z^n=\frac{1}{\sqrt{1-z}}, \quad |z|<1,
\]
implies that
\[
\frac{1}{2 \pi i} \int_{-i\infty}^{i\infty} \left( \frac{1}{2} \right)_{\!s} \Gamma(-s) (-z)^s ds=\frac{1}{\sqrt{1-z}}, \quad z \in \mathbb{C}-[1,\infty).
\]
Hence, we have
\[
\frac{1}{2 \pi i} \int_{-i\infty}^{i\infty} \frac{\left( \frac{1}{2} \right)_{\!s} \left( \frac{1}{4}+\frac{3t}{2} \right)_{\!s} \left( \frac{3}{4}+\frac{3t}{2} \right)_{\!s}}{(1+t)_{\!s} (1+2t)_{\!s}}  \frac{\left( \frac{1}{6} \right)_t \left( \frac{5}{6} \right)_t}{(1)_t^2} (5s+6t+1) \Gamma(-s)\left( \frac{4}{3} \right)^{2s}ds=\frac{\sqrt{3}}{\pi},
\]
or equivalently
\[
\frac{1}{2 \pi i} \int_{-i\infty}^{i\infty} \frac{\left( \frac{1}{2} \right)_{\!s} \left( \frac{1}{4}+\frac{3t}{2} \right)_{\!s} \left( \frac{3}{4}+\frac{3t}{2} \right)_{\!s}}{(1+t)_{\!s} (1+2t)_{\!s}}   (5s+6t+1) \Gamma(-s)\left( \frac{4}{3} \right)^{2s}ds=\frac{\sqrt{3}}{\pi} \frac{(1)_t^2}{\left( \frac{1}{6} \right)_t \left( \frac{5}{6} \right)_t}.
\]
Finally, substituting $t=0$, we see that
\begin{equation}\label{for5s1}
\frac{1}{2 \pi i} \int_{-i\infty}^{i\infty} \frac{\left( \frac{1}{2} \right)_{\!s} \left( \frac{1}{4} \right)_{\!s} \left( \frac{3}{4} \right)_{\!s}}{(1)_{\!s}^2}(5s+1) \Gamma(-s)\left( \frac{4}{3} \right)^{2s}ds=\frac{\sqrt{3}}{\pi}.
\end{equation}
\par It is very convenient to write the Barnes integral in hypergeometric notation. By the definition of hypergeometric series, we see that for $-1 \leq z < 1$, we have
\[
\sum_{n=0}^{\infty} \frac{\left( \frac{1}{2} \right)_n \left( s \right)_n \left( 1-s \right)_n}{(1)_n^3}z^n =
{}_3F_2\biggl(\begin{matrix}
\frac12, & s, & 1-s \\
& 1, & 1 \end{matrix} \biggm| z\biggr)
\]
and
\[
\sum_{n=0}^{\infty} \frac{\left( \frac{1}{2} \right)_n \left( s \right)_n \left( 1-s \right)_n}{(1)_n^3} nz^n = \frac12 s(1-s)z \,
{}_3F_2\biggl(\begin{matrix}
\frac32, & 1+s, & 2-s \\
& 2, & 2 \end{matrix} \biggm| z\biggr),
\]
where the notation on the right side stands for the analytic continuation of the series on the left. Hence, we can write (\ref{for5s1}) in the form
\[
{}_3F_2\biggl(\begin{matrix}
\frac12, & \frac14, & \frac34 \\
& 1, & 1 \end{matrix} \biggm| \frac{-16}{9} \biggr)-\frac{5}{6} \,
{}_3F_2\biggl(\begin{matrix}
\frac32, & \frac54, & \frac74 \\
& 2, & 2 \end{matrix} \biggm| \frac{-16}{9} \biggr)=\frac{\sqrt{3}}{\pi}.
\]
\par If, instead of integrating to the right side, we integrate (\ref{for5s1}) along a contour which closes the path $\mathcal{P}$ with a semicircle of center $s=0$ taken to the left side with an infinite radius, then we have poles at $s=-n-1/2$, at $s=-n-1/4$ and at $s=-n-3/4$ for $n=0,1,2, \dots$, and we obtain
\begin{multline}\nonumber
\frac{\sqrt{3}}{2} \sum_{n=0}^{\infty} \frac{\left( \frac12\right)_n^3}{(1)_n \left( \frac34\right)_n
\left( \frac54\right)_n}(10n+3)(-1)^n \left( \frac{3}{4} \right)^{2n}
\\ \nonumber - \frac{\sqrt{2} \: \pi^2}{8 \; \Gamma \! \left( \frac34 \right)^4}
\sum_{n=0}^{\infty} \frac{\left( \frac14\right)_n^3}{(1)_n \left( \frac12 \right)_n \left( \frac34\right)_n}(20n+1) (-1)^n \left( \frac{3}{4} \right)^{2n}
\\ \nonumber - \frac{3\sqrt{2} \: \Gamma \! \left( \frac34 \right)^4}{16 \: \pi^2}
\sum_{n=0}^{\infty} \frac{\left( \frac34\right)_n^3}{(1)_n
\left( \frac32 \right)_n\left( \frac54\right)_n}(20n+11)(-1)^n \left( \frac{3}{4} \right)^{2n}=1.
\end{multline}
which is an identity relating three convergent series.

\section{Other examples}
In a similar way we can prove other identities of the same kind, for example,
\[
\frac{1}{2 \pi i} \int_{-i\infty}^{i\infty} \frac{\left( \frac{1}{2}+t \right)_{\!s}^3 \left( \frac{1}{2} \right)_{\!s}^2 }{(1+t)_{\!s}^3(1+2t)_{\!s}}(10s^2+6s+1+14st+4t^2+4t) \Gamma(-s)2^{2s}ds=\frac{4}{\pi^2} \frac{(1)_t^4}{\left( \frac{1}{2} \right)_t^4},
\]
\[
\frac{1}{2 \pi i} \int_{-i\infty}^{i\infty} \frac{\left( \frac{1}{2} \right)_{\!s} \left( \frac{1}{2}+t \right)_{\!s}^2 }{(1)_{\!s}(1+2t)_{\!s}}(3s+2t+1) \Gamma(-s) 2^{3s}ds=\frac{1}{\pi} \frac{(1)_t}{\left( \frac{1}{2} \right)_t},
\]
and
\begin{multline}\nonumber
\frac{1}{2 \pi i} \int_{-i\infty}^{i\infty}  \frac{\left( \frac{1}{2} \right)_{\!s} \left( \frac{1}{2}+2t \right)_{\!s} \left( \frac{1}{3}+t \right)_{\!s} \left( \frac{2}{3}+t \right)_{\!s}}{\left(\frac{1}{2}+\frac{t}{2}\right)_{\!s}  \left(1+\frac{t}{2}\right)_{\!s} (1+t)_{\!s} } \\ \times \frac{(15s+4)(2s+1)+t(33s+16)}{2s+t+1} \Gamma(-s)2^{2s}ds  = \frac{3\sqrt{3}}{\pi} \frac{1}{2^{6t}} \frac{(1)_t^2}{\left( \frac{1}{4} \right)_t \left( \frac{3}{4} \right)_t}.
\end{multline}
In the two last examples the hypothesis of Weierstrass theorem fail and hence we cannot apply it, but we obtain the sum using Meurman's periodic version of Carlson's theorem \cite[p. 39]{bailey} which asserts that if $H(z)$ is a periodic entire function of period $1$ and there is a real number $c<2\pi$ such that $H(z)=\mathcal{O}(\exp(c |Im(z)|))$ for all $z \in \mathbb{C}$, then $H(z)$ is constant \cite[Appendix]{almkvist} and \cite[Thm. 2.3]{guillera-morehypiden}. In the second and third examples we determine the constants $1/\pi$ and $3 \sqrt{3} / \pi$ taking $t=1/2$ and $t=-1/3$ respectively. Substituting $t=0$ in the above examples, we obtain respectively
\begin{equation}\label{ej1}
\frac{1}{2 \pi i} \int_{-i\infty}^{i\infty} \frac{\left( \frac{1}{2} \right)_{\!s}^5}{(1)_{\!s}^4}(10s^2+6s+1) \Gamma(-s)2^{2s}ds=\frac{4}{\pi^2},
\end{equation}
\begin{equation}\label{ej2}
\frac{1}{2 \pi i} \int_{-i\infty}^{i\infty} \frac{\left( \frac{1}{2} \right)_{\!s}^3}{(1)_{\!s}^2}(3s+1) \Gamma(-s) 2^{3s}ds=\frac{1}{\pi},
\end{equation}
and
\begin{equation}\label{ej3}
\frac{1}{2 \pi i} \int_{-i\infty}^{i\infty}  \frac{\left( \frac{1}{2} \right)_{\!s} \left( \frac{1}{3} \right)_{\!s} \left( \frac{2}{3} \right)_{\!s}}{\left(1\right)_{\!s}^2} (15s+4) \Gamma(-s)2^{2s}ds=\frac{3\sqrt{3}}{\pi}.
\end{equation}
Using hypergeometric notation, we can write (\ref{ej1}), (\ref{ej2}) and (\ref{ej3}) respectively in the following forms:
\[
{}_5F_4\biggl(\begin{matrix}
\frac12, & \frac12, & \frac12, & \frac12, & \frac12 \\
& 1, & 1, & 1, & 1 \end{matrix} \biggm| -4 \biggr)
-\frac{3}{4} \, {}_5F_4\biggl(\begin{matrix}
\frac32, & \frac32, & \frac32 & \frac32 & \frac32 \\
& 2, & 2, & 2, & 2  \end{matrix} \biggm| -4 \biggr)
\]
\[
\hskip 6.5cm -\frac{5}{4} \, {}_5F_4\biggl(\begin{matrix}
\frac32, & \frac32, & \frac32 & \frac32 & \frac32 \\
& 2, & 2, & 2, & 1  \end{matrix} \biggm| -4 \biggr)=\frac{4}{\pi^2},
\]
\[
{}_3F_2\biggl(\begin{matrix}
\frac12, & \frac12, & \frac12 \\
& 1, & 1 \end{matrix} \biggm| -8 \biggr)-
3 \, {}_3F_2\biggl(\begin{matrix}
\frac32, & \frac32, & \frac32 \\
& 2, & 2 \end{matrix} \biggm| -8 \biggr)=\frac{1}{\pi},
\]
and
\[
4 \, {}_3F_2\biggl(\begin{matrix}
\frac12, & \frac13, & \frac23 \\
& 1, & 1 \end{matrix} \biggm| -4 \biggr)-\frac{20}{3} \,
{}_3F_2\biggl(\begin{matrix}
\frac32, & \frac43, & \frac53 \\
& 2, & 2 \end{matrix} \biggm| -4 \biggr)=\frac{3\sqrt{3}}{\pi}.
\]
Related applications of the WZ-method for Barnes-type integrals are for example in \cite[Sect. 5.2]{bailey-ising} and \cite{stan}.

\section{The dual of a ``divergent'' Ramanujan-type series}

The WZ duality technique \cite[Ch. 7]{petkovsek} allows to transform pairs which lead to divergences into pairs which lead to convergent series. To get the dual $\widehat{G}(n,k)$ of $G(-n,-k)$, we make the following changes:
\[
(a)_{-n} \rightarrow \frac{(-1)^n}{(1-a)_n}, \quad (1)_{-n} \rightarrow \frac{n(-1)^n}{(1)_n}, \quad (a)_{-k} \rightarrow \frac{(-1)^k}{(1-a)_k}, \quad (1)_{-k} \rightarrow \frac{k(-1)^k}{(1)_k}.
\]

\subsection*{Example 1}
The package EKHAD certifies the pair
\begin{equation}
F(n,k)=U(n,k) \frac{2n^2}{2n+k}, \qquad G(n,k)=U(n,k) \frac{6n^2+2n+k+4nk}{2n+k},
\end{equation}
where
\[
U(n,k)=\frac{\left(\frac{1}{2} \right)_n^2 \left(1+\frac{k}{2}\right)_n \left(\frac{1}{2}+\frac{k}{2}\right)_n}{(1)_n^2(1+k)_n^2}
\frac{\left(\frac{1}{2}\right)_k}{(1)_k} \, \, 4^n= \frac{(2n)!^2(2n+k)!(2k)!}{n!^4 k! (n+k)!^2} \frac{1}{16^n 4^k}.
\]
We cannot use this WZ-pair to obtain a Ramanujan-like evaluation because, as $z>1$, the corresponding series and also the corresponding Barnes integral are both divergent. However, we will see how to use it to evaluate a related convergent series. What we will do is to apply the WZ duality technique. Thus, if we take the dual of $G(-n,-k)$ and replace $k$ with $k-1$, we obtain
\[
\widehat{G}(n,k)=\frac{1}{U(n,k)} \frac{2(2k-1)(2n+k)}{n^2(n+k)^2(n+k-1)^2} (6n^2-6n+1-k+4nk),
\]
and EHKAD finds its companion
\[ \widehat{F}(n,k)=\frac{1}{U(n,k)} \frac{-2(2n+k)(2n+k-1)(2n-1)^2}{n^2(n+k)^2(n+k-1)^2}. \]
Applying Zeilberger's formula
\[
\sum_{n=j}^{\infty} (\widehat{F}(n+1,n)+\widehat{G}(n,n))=\sum_{n=j}^{\infty} \widehat{G}(n,j)
\]
with $j=1$, we obtain
\begin{equation}\label{identidad}
\sum_{n=1}^{\infty} \left( \frac{16}{27} \right)^n \frac{(1)_n^3}{\left(\frac{1}{2} \right)_n \left( \frac{1}{3} \right)_n \left(\frac{2}{3} \right)_n} \frac{11n-3}{n^3}=16 \sum_{n=1}^{\infty} \frac{1}{4^n} \frac{(1)_n^3}{\left(\frac{1}{2} \right)_n^3} \frac{3n-1}{n^3}.
\end{equation}
The series in (\ref{identidad}) are dual to Ramanujan-type ``divergent'' series, and in \cite[p. 221]{gui-hyperiden} we proved that the series on the right side is equal to $\pi^2/2$. Hence
\begin{equation}\label{zhi}
\sum_{n=1}^{\infty} \left( \frac{16}{27} \right)^n \frac{(1)_n^3}{\left(\frac{1}{2} \right)_n \left( \frac{1}{3} \right)_n \left(\frac{2}{3} \right)_n} \frac{11n-3}{n^3}=8 \pi^2.
\end{equation}
Formula (\ref{zhi}), as well as other similar formulas, was conjectured in \cite[Conj 1.4]{sun} by Zhi-Wei Sun.

\subsection*{Example 2}

The package EKHAD certifies the pair
\begin{align}
F(n,k) & =U(n,k) \frac{64n^3}{(2k+1)(2n-2k+1)}, \nonumber \\
G(n,k) & =U(n,k) \frac{(2n+1)^2(11n+3)-12k(2n^2+3nk+n+k)}{(2n+1)^2}, \nonumber
\end{align}
where
\[
U(n,k)=\frac{\left(\frac{1}{2}-k \right)_n \left(\frac12+k\right)_n^2 \left(\frac{1}{3}\right)_n \left(\frac{1}{3}\right)_n}{(1)_n^3 \left( \frac12 \right)_n^2} \, \, \left( \frac{27}{16} \right)^n.
\]
Taking the dual $\widehat{G}(n,k)$ of $G(-n,-k)$, replacing $n$ with $n+x$ and applying Zeilberger's theorem
\[
\sum_{n=0}^{\infty} \widehat{G}(n+x,0)=\lim_{k \to \infty} \sum_{n=0}^{\infty}\widehat{G}(n+x,k)+\sum_{k=0}^{\infty} \widehat{F}(x,k), \]
where $\widehat{F}(n,k)$ is the companion of $\widehat{G}(n,k)$, we obtain
\begin{multline}\nonumber
\sum_{n=0}^{\infty} \frac{(1+x)_n^3}{\left( \frac12+x\right)_n\left( \frac13+x\right)_n\left( \frac23+x\right)_n}\left( \frac{16}{27} \right)^n \frac{11(n+x)-3}{(n+x)^3} \\= \frac{6(3x-1)(3x-2)}{x^3(2x-1)}\sum_{k=0}^{\infty} \frac{\left( \frac12 \right)_k \left( \frac32-x \right)_k}{\left( \frac12+x \right)_k^2}.
\end{multline}
Taking $x=1$ we again obtain (\ref{zhi}).

\section*{Acknowledgment}

I thank W. Zudilin for encouraging me to use the Barnes integral representation to obtain the sum of ``divergent'' Ramanujan-type series. I thank G. Almkvist for suggesting to me to integrate on the left side to get identities with convergent series. I am also grateful to Jonathan Sondow for several helpful comments.

\vskip 0.35cm

\end{document}